\documentclass{amsart}
\usepackage{geometry} 
\usepackage{url}
\usepackage{amsthm}
\usepackage{amssymb}
\usepackage{latexsym}
\usepackage{amsfonts}
\usepackage{amsmath}
\usepackage{amstext}
\usepackage{accents}
\usepackage{cancel}
\usepackage{eucal}
\usepackage{amscd}
\usepackage{hyperref}
\usepackage[shortalphabetic]{amsrefs}
\usepackage{tikz}

\begin{document}
\title[The Length-Constrained Ideal Curve Flow]{\textbf{The Length-Constrained Ideal Curve Flow}}

\author[J. A. McCoy]{James A. McCoy}
\address{Priority Research Centre Computer Assisted Research Mathematics and Applications, School of Mathematical and Physical Sciences, University of Newcastle, University Drive, Callaghan, NSW 2308, Australia}
\email{James.McCoy@newcastle.edu.au}

\author[G. E. Wheeler]{Glen E. Wheeler}
\address{Institute for Mathematics and its Applications, School of Mathematics and Applied Statistics, University of Wollongong, Northfields Avenue, Wollongong, NSW 2522, Australia}
\email{glenw@uow.edu.au}

\author[Y. Wu]{Yuhan Wu}
\address{Institute for Mathematics and its Applications, School of Mathematics and Applied Statistics, University of Wollongong, Northfields Avenue, Wollongong, NSW 2522, Australia}
\email{yw120@uowmail.edu.au}

\subjclass[2000]{53E40, 58J35}

\begin{abstract}
A recent article \cite{AMWW} considered the so-called `ideal curve flow', a sixth order curvature flow that seeks to deform closed planar curves to curves with least variation of total geodesic curvature in the $L^2$ sense.  Critical in the analysis there was a length bound on the evolving curves.  It is natural to suspect therefore that the length-constrained ideal curve flow should permit a more straightforward analysis, at least in the case of small initial `energy'.  In this article we show this is indeed the case, with suitable initial data providing a flow that exists for all time and converges smoothly and exponentially to a multiply-covered round circle of the same length and winding number as the initial curve.
\end{abstract}

\keywords{geometric evolution equation, higher order parabolic partial differential equation, higher order curvature flow}

 \thanks{The research of the first and second author was partially supported by Discovery Grant DP180100431 of the Australian Research Council.  The research of the third author was supported by a University of Wollongong Faculty of Engineering and Information Sciences Postgraduate research scholarship.  The authors are grateful for these sources of support.}

\maketitle

\section{Introduction}
\newtheorem{main}{Theorem}[section]

In the recent article \cite{AMWW}, the first and second authors, together with Andrews and V-M Wheeler, considered, for smooth, closed, planar curves $\gamma: \mathbb{S} \rightarrow \mathbb{R}^2$ the $L^2$ gradient flow of the energy
\begin{equation} \label{E:energy}
  E\left[ \gamma\right] = \frac{1}{2} \int_\gamma k_s^2 ds \mbox{,}
  \end{equation}
where $k_s$ denotes the first arclength derivative of the curvature $k$ of $\gamma$.  This energy appears, for example, in computer aided design \cite{HT}.  Writing now $\gamma: \mathbb{S} \times \left[ 0, T\right) \rightarrow \mathbb{R}^2$, the corresponding $L^2$ gradient flow is given explicitly by the evolution equation
\begin{equation} \label{E:ideal}
\partial_t \gamma = \left( k_{s^4} + k^2 k_{ss} - \frac{1}{2} k\, k_s^2 \right) \nu =: G\left[ k \right] \nu
\end{equation}
with given initial data 
$$\gamma\left( \cdot, 0 \right) = \gamma_0 \mbox{.}$$
Above and in general we will use the notation $k_{s^m}$ to denote the $m$-th iterated derivative of curvature with respect to arclength $s$.	

Multiply-covered circles are obviously absolute minimisers of \eqref{E:energy} but this property is independent of their radius, so understanding the behaviour of the length functional $L\left[ \gamma\right]$ under the flow \eqref{E:ideal} is essential.  In \cite{AMWW} it was shown that a solution to \eqref{E:ideal} always exists for all time and if the length is uniformly bounded then the solution converges exponentially fast in the $C^\infty$-topology to a standard, round $\omega$-circle, where 
\begin{equation} \label{E:wn}
  \omega= \frac{1}{2\pi} \int_{\gamma_0} k_0\, ds_0
\end{equation}
is the winding number of $\gamma_0$.  Moreover it was shown that if the initial scale-invariant energy $\left( L^3 E\right) \left[ \gamma_0 \right]$ is sufficiently small, then $L\left[ \gamma_t\right]$ is uniformly bounded under the flow \eqref{E:ideal} giving one criterion for which the convergence statement holds.

In view of the above, it is natural to expect that adding a constraint to the flow speed that ensures the length of the evolving curve remains constant could simplify the analysis, at least in the case of small `energy' (in this article we will continue to refer to \eqref{E:energy} as a corresponding energy to our constrained flow, even though our flow is no longer the $L^2$-gradient flow of the energy and monotonicity of the quantity is not clear in general).  We show that indeed the length constraint significantly simplifies the analysis to provide exponential convergence of initial small `energy' curves to round $\omega$-circles.  Specifically we prove the following result:

\begin{main} \label{T:main}
Let $\gamma_0: \mathbb{S}\to \mathbb{R}^2$ be a regular smooth immersed closed curve with length $L_0$ and winding number $\omega$, satisfying the condition
$$E[\gamma_0]\leq\varepsilon \mbox{,}$$
where $\varepsilon$ depends only on $L_0$ and $\omega$.  Then with initial curve $\gamma\left( \cdot, 0\right) =\gamma_0$, the \emph{length-constrained ideal curve flow}
\begin{equation} \label{r1}
\partial_t^{\perp}\gamma=k_{s^4}+k^2k_{ss}-\frac{1}{2}kk_s^2+h(t)
\end{equation}
where
\begin{equation}\label{h2}
h(t)=\frac{1}{2\pi w}\left(-\int k_{ss}^2ds+\frac{7}{2}\int k_s^2k^2ds\right)
\end{equation}
has a smooth solution, unique up to parametrisation, that exists for all time and converges exponentially in the $C^\infty$ topology to a multiply-covered $\omega$-circle of radius $\frac{L_0}{2\pi\omega}$.  
\end{main}

\noindent \emph{Remarks:}
\begin{enumerate}
  \item Above $\partial_t^\perp \gamma$ denotes the normal component of the time derivative of $\gamma$.  The normal direction is chosen such that the flow \eqref{r1} is parabolic in the generalised sense.
  \item Keeping the definition of $G$ from \eqref{E:ideal}, we may rewrite \eqref{r1} and \eqref{h2} respectively as
\begin{equation} \label{E:speed}
  \partial_t^{\perp}\gamma=G+h(t) 
  \end{equation}
  and
  $$h(t)=-\frac{1}{2\pi \omega} \int k\, G \, ds \mbox{.}$$
  Each of these will be useful later.\\
  \end{enumerate}

The structure of this article is as follows.  We will complete this introductory section with a brief review of relevant previous work.  In Section \ref{S:prelim} we provide the tools that will be required in our analysis.  In Section \ref{S:evlneqns} we give the relevant evolution equations and some preliminary results.  In Section \ref{S:LTE} we complete the proof of Theorem \ref{T:main}

Of higher order curvature flows of closed planar curves, the fourth order cases of the curve diffusion flow and the elastic flow are the most well-studied.  Important early work includes \cite{K96, LS84, LS85, L89, L98, P, W93, W95}.  There is also the work of Giga and Ito \cite{GI1, GI2} showing that curves may lose embeddedness or convexity under the curve diffusion flow.  In 2002, Dzuik, Kuwert and Sch\"{a}tzle studied the elastic flow of curves (in arbitrary codimension), establishing some results including interpolation inequalities that are applicable broadly and fundamental here.  Later the curve diffusion flow of closed planar curves was studied by the second author where it was shown that initial curves suitably close to a circle in the sense of the normalised $L^2$-norm of the oscillation of curvature converge exponentially to a circle \cite{W1}.  Okabe considered an area-preserving elastic flow in \cite{O07} and with uniform high pressure acting on the curve in \cite{O08}.
Okabe and Novaga considered evolving curves of infinite length in \cite{NO14}.  The curve diffusion flow has the basic property that the signed enclosed area is preserved under the flow while the length of the evolving  curve does not increase.  More recently the authors of this article considered the length constrained curve diffusion flow \cite{MWW} showing an analogous result.  Under that flow, the length of the evolving curve is fixed and the signed enclosed area does not decrease.  Whereas the global term in that flow has a clear sign, the term \eqref{h2} for the flow considered in this article does not.  

Sixth and higher order flows of closed planar curves have been less well-studied, despite applications in imaging and computer design \cite{HT, UW}.  As mentioned earlier, the main counterpart for this article is \cite{AMWW}, where the sixth order flow \eqref{E:ideal} was considered.  The polyharmonic flow of closed plane curves was considered in \cite{PW16} and an anisotropic analogue in \cite{PW19}.  If we move to the setting of curves with (generalised) Neumann boundary conditions on parallel lines, more is known, often due to the simplification of the Poincar\'{e}-Sobolev-Wirthinger inequalities in this setting.  The interpolation inequality of \cite{DKS02} was generalised to the case of
evolving curves with boundary by Dall'Acqua and Pozzi \cite{DP14}; in that article the
authors focus on the so-called Willmore-Helfrich flow of curves with natural
boundary conditions.  The second author and V-M Wheeler considered the fourth order cases of curve diffusion and the elastic flow in \cite{WW}.  In this setting, the analog of \cite{AMWW} was also considered by the authors \cite{MWW1} and this was further generalised in \cite{MWW2} to both polyharmonic flow with small initial curvature in $L^2$ and small-energy $L^2$ gradient flows of $\int k_{s^m}^2 ds$.

\section{Preliminaries} \label{S:prelim}
\newtheorem{PSW}{Lemma}[section]
\newtheorem{DKS}[PSW]{Lemma}
\newtheorem{ksup}[PSW]{Lemma}
\newtheorem{G2est}[PSW]{Lemma}
\newtheorem{G2new}[PSW]{Corollary}
Let $\gamma_0: \mathbb{R}\to \mathbb{R}^2$ be a smooth embedded (or immersed) regular curve. We say $\gamma$ is periodic with period $P$ if there exists a vector $V\in \mathbb{R}^2$ and a positive number $P$ such that, for all $l\in \mathbb{N}$ and every $u\in \mathbb{R}$
$$\gamma(u+P)=\gamma(u)+V \mbox{ and } \partial_u^\ell\gamma(u+P)=\partial_u^\ell\gamma(u) \mbox{.}$$
Here $\partial_u^\ell$ denotes the $\ell$th iterated derivative of $\gamma$. If $V=0$ then $\gamma$ is closed and we may rewrite $\gamma: \mathbb{S}\to \mathbb{R}^2$. The length of $\gamma$ is 
$$L[\gamma] :=\int_0^P\vert \gamma'(u)\vert du \mbox{.}$$

Throughout this article we will keep our evolving curves $\gamma$ parametrised by arc length $s$, defined via
$$s\left( u \right) = \int_{u_0}^u \left| \gamma'\left( u\right) \right| du \mbox{,}$$
for a choice of $u_0 \in \mathbb{S}$.

In our analysis we will frequently use the following Poincar\'{e}-Sobolev-Wirtinger [PSW] inequalities. For proofs of these see for example Appendix A of \cite{PW16}. 
 
\begin{PSW} \label{T:PSW}
Suppose $f: \left[ 0, P\right] \to \mathbb{R}$ is a absolutely continuous.  Then if $\int _0^P f(x)dx=0$ we have 
\begin{enumerate}
  \item[\textnormal{(i)}]
$$\int _0^P f^2(x)dx\leq \frac{P^2}{4\pi^2}\int _0^P  f'(x)^2dx$$
with equality if and only if $f(x)=a \sin \left(\frac{2\pi x}{P}+b\right)$ for arbitrary constants $a$ and $b$; 
  \item[\textnormal{(ii)}]
  $$\Vert f\Vert_\infty^2\leq \frac{P}{2\pi}\int_0^P f'(x)^2dx$$
  \end{enumerate}
\end{PSW}
We also need some interpolation inequalities from \cite{DKS02}. We first set up some notation. For normal tensor fields $S$ and $T$ we denote by $S\star T$ any linear combination of $S$ and $T$. In our setting, $S$ and $T$ will be simply curvature $k$ or its arclength derivatives. We denote by $P_n^m$ any linear combination of terms of type $\partial_s^{i_1}l\star \partial_s^{i_2}l\star ... \star \partial_s^{i_n}l$, where $m=i_1+...+i_n$ is the total number of derivatives.

The following interpolation inequality for closed curves appears in \cite{DKS02}. 

\begin{DKS}\label{i13}
Let $\gamma:\mathbb{S} \to \mathbb{R}^2$ be a smooth closed curve. Then for any term $P_n^m(k)$ with $n\geq 2$ that contains derivatives of $k$ of order at most $l-1$,
$$\int _\gamma\vert P_n^m(k)\vert ds\leq cL^{1-m-n}\Vert k\Vert_2^{n-p}\Vert k\Vert_{l,2}^p$$
where $p=\frac{1}{l}\left(m+\frac{1}{2}n-1\right)$ and $c=c(l,m,n)$. Moreover, if $m+\frac{n}{2}<2l+1$ then $p<2$ and for any $\varepsilon>0$,
$$\int_\gamma\vert P_n^m(k)\vert ds\leq \varepsilon\int_\gamma \vert \partial_{s^\ell}k\vert^2ds+c\varepsilon^{\frac{-p}{2-p}}\left(\int_\gamma k^2ds\right)^{\frac{n-p}{2-p}}+c\left(\int_\gamma k^2ds\right)^{m+n-1} \mbox{.}$$
\end{DKS}

In the above, $\Vert \cdot\Vert_s$ and $\Vert \cdot \Vert_{m,2}$ denote scale-invariant norms, for example
$$\Vert k\Vert_2=\Vert k\Vert_{0,2}=L^{\frac{1}{2}}\left(\int k^2ds\right)^{\frac{1}{2}}$$
and 
$$\Vert k\Vert_{1,2}=L^{\frac{1}{2}}\left( \int k^2ds\right)^{\frac{1}{2}}+L^{\frac{3}{2}}\left(\int k_s^2ds\right)^{\frac{1}{2}}.$$
Under the flow in this article $L=L_0$ is constant, so the $L$ factors above are harmless.  In our estimates we will also allow the constants $c$ to vary from line to line where they depend only on absolute quantities like $n$, $m$ and, for this flow $L_0$. Of course when $L_0$ is contained within our constants it is no longer possible to track scaling through the estimates.  In our setting, since $L_0$ is constant, so is the average curvature under the flow:
$$\overline{k} := \frac{1}{L} \int_\gamma k\, ds = \frac{2\pi \omega}{L_0} \mbox{.}$$
(In Lemma \ref{T:wn} we show that $\omega$ is also constant under the flow \eqref{r1}.)\\

The next elementary inequality we state is \cite[Lemma 6.2]{AMWW}.

\begin{ksup} \label{T:ksup}
For any immersed curve $\gamma: \mathbb{S} \rightarrow \mathbb{R}^2$ with winding number $\omega$,
$$L \left\| k \right\|_\infty \leq \sqrt{L^3 \left\| k_s \right\|_2^2} + 2\omega \pi \mbox{.}$$
\end{ksup}

\noindent \textbf{Remarks:} 
\begin{enumerate}
  \item In our case since $L=L_0$ is constant we see that $\left\| k \right\|_{\infty}$ is bounded as soon as we have control on $\left\| k_s \right\|_2^2$.\\
  \item Another obvious corollary of Lemma \ref{T:ksup}, given $L_0$ is constant, is
\begin{equation*} 
  \left\| k+ \overline{k} \right\|_\infty \leq \left\| k \right\|_\infty + \overline{k} \leq \sqrt{ 2 L_0 E} + \frac{4\pi}{L_0} \mbox{.}
  \end{equation*}
\end{enumerate}

We also have from \cite[Proposition 7.1]{AMWW} the following crucial flow-independent estimate, that may be proved using a Fourier series argument.

\begin{G2est} \label{T:G2est}
For a smooth, closed curve of length $L$ and winding number $\omega$, there exist positive constants $c_\omega$ depending only on $\omega$ and $c$, depending on $\omega$ and $L$, such that
$$\int \left( k_{s^4}^2 + \overline{k}^2 k_{ss}\right)^2 ds \geq c_\omega \int k_{s^4}^2 ds - c\left( \omega, L\right) \left( \int k_s^2 ds\right)^2 \mbox{.}$$
\end{G2est}

In this article we make use of the following consequence, again flow-independent, where the constants are different.

\begin{G2new} \label{T:G2new}
For a smooth, closed curve of length $L$ and winding number $\omega$, there exist positive constants $c_\omega$ depending only on $\omega$ and $c$, depending on $\omega$ and $L$, such that
$$\int G^2 ds \geq c_\omega \int k_{s^4}^2 ds - c\left( \omega, L\right) E^2 \left( 1 + E^3 \right)\mbox{.}$$
\end{G2new}

\noindent \textbf{Proof:} The proof follows that of \cite[Proposition 7.2]{AMWW}.  Specifically, using Lemma \ref{T:G2est}, Lemma \ref{T:ksup} and the Gagliardo-Nirenberg Sobolev inequality exactly as there, we obtain, for positive constants $c_1$ and $c_2$ depending only on $\omega$,
$$\int G^2 ds \geq  c_1 \int k_{s^4}^2 ds - c_2 \left\{ L^{-3} E^2 + \left[ \left( L^3 E\right)^{\frac{3}{2}} + \omega^3\right] \left( L^{-\frac{3}{2}} E^{\frac{5}{2}} \right) + \left[ \left( L^3 E\right)^2 + \omega^4 \right] E^3 \right\} \mbox{.}$$
The result follows by expanding and estimating.\hspace*{\fill}$\Box$\\

\noindent \emph{Remark:} Because $L$ is constant along the flow \eqref{r1}, the form of inequality in Corollary \ref{T:G2new} is sufficient for our purposes.

\section{Evolution equations and first estimates} \label{S:evlneqns}
\newtheorem{commute}{Lemma}[section]
\newtheorem{evlneqns}[commute]{Lemma}
\newtheorem{wn}[commute]{Lemma}
\newtheorem{Lconst}[commute]{Lemma}
\newtheorem{Lbound}[commute]{Corollary}
\newtheorem{Fbar}[commute]{Lemma}
\newtheorem{kslevln}[commute]{Lemma}

Let us first collect some fundamental evolution equations that will be needed in our analysis.  Derivations of these for general normal flow speed $F$ are straightforward, as in \cite{W1, WW}, for example.  Since the arclength of the evolving curve is changing in time, we first need the commutator of arclength and time derivatives.

\begin{commute} \label{p7}
Under flow by normal speed $F$, the commutator of arclength and time derivatives is given via
$$\partial_t \partial_s = \partial_s \partial_t + k F \,\partial_s$$
and the measure $ds$ evolves according to
$$\partial_t ds = - k F\, ds \mbox{.}$$
\end{commute}

Next we give  evolution equations for length and for arclength curvature derivatives under flow by normal speed $F$:

\begin{evlneqns}\label{d12}
\begin{enumerate}
  \item[\textnormal{(i)}] $\frac{d}{d t}L=-\int_\gamma kFds=0$;
\item[\textnormal{(ii)}] For each $l=0,1,2, ..., $
$$\frac{\partial}{\partial t}k_{s^l}=F_{s^{l+2}}+\sum_{j=0}^l \partial_{s^j}\left(kk_{s^{l-j}}F\right) \mbox{;}$$
\end{enumerate}
\end{evlneqns}

It is a quite general property for evolving curves that under various conditions the winding number of $\gamma_t$ remains constant.  We provide the proof for our setting.

\begin{wn} \label{T:wn}
Under flow by normal speed $F$, the winding number of closed curve $\gamma$ remains constant.
\end{wn}

\noindent \textbf{Proof:}. As $\omega$ is given via \eqref{E:wn}, we compute using Lemma \ref{p7} and Lemma \ref{d12} (ii), 
$$\frac{d}{dt}\int_\gamma k\, ds=\int_\gamma \frac{\partial}{\partial t}k\, ds+\int_\gamma k\cdot (-kF)\, ds
=\int_\gamma \left(F_{ss}+k^2F-k^2F\right)ds=0$$
The result follows. \hspace*{\fill}$\Box$\\

Short-time existence of a unique solution modulo reparametrisation to the flow \eqref{r1} follows by standard theory, writing the solution as a graph over the initial curve for a short time.  We refer the reader for a detailed discussion in a more general setting in \cite{B}, noting that the term $h\left( t\right)$ can be handled by a standard fixed point argument similarly as in \cite{M}, for example.  We refer also to \cite{W} by the third author of short time existence in several settings.

In fact, a standard argument shows that the solution to \eqref{r1} exists so long as $\int k^2 ds$ is finite.  The argument uses Lemma \ref{i13} to bound all $L^2$ norms of the curvature derivatives in terms of $\int k^2 ds$, so if this quantity is finite then the solution can be extended via short-time existence.  A similar argument is used in Section \ref{S:LTE} of the present article, although there we use also the decay of $E\left[ \gamma_t\right]$ which requires a smallness condition.

It is straightforward to verify that the choice of $h\left( t\right)$ given by \eqref{h2} ensures the length $L\left[ \gamma_t\right]$ is constant along the flow.

\begin{Lconst}
Under the flow \eqref{r1}, while a solution exists, it satisfies
$$L\left[ \gamma_t\right] = L_0 \mbox{.}$$
\end{Lconst}

\noindent \textbf{Proof:} We compute
\begin{multline*}
 \frac{d}{dt} L \left[ \gamma_t\right] = -\int_\gamma k F \, ds = -\int k\, k_{s^4}ds - \int k^3 k_{ss} ds + \frac{1}{2} \int k^2 k_s^2 ds + h\left( t\right) \int k\, ds \\
 = - \int k_{ss}^2 ds + \frac{7}{2} \int k^2 k_s^2 ds - 2\pi\omega h\left( t\right) \mbox{;}
 \end{multline*}
 the result follows in view of \eqref{h2}.\hspace*{\fill}$\Box$\\

We complete this section with the evolution for arclength curvature derivatives in $L^2$ under the flow \eqref{r1}.  This is of course similar to that in \cite{AMWW} except here we have an additional term arising from $h\left( t\right)$.

\begin{kslevln} \label{T:kslevln}
Under the flow \eqref{r1},
$$\frac{d}{dt} \int k_{s^\ell}^2 ds = - 2 \int k_{s^{\ell+3}}^2 ds + 2 \int k_{s^\ell} \left[ P_3^{\ell+4}\left( k\right) + P_5^{\ell+2}\left( k\right) \right] ds - h\left( t\right) \int k_{s^\ell}^2 k\, ds \mbox{.}$$
\end{kslevln}

\section{The estimates} \label{S:LTE}
\newtheorem{Emonotone}{Proposition}[section]
\newtheorem{Edecay}[Emonotone]{Corollary}
\newtheorem{kbound}[Emonotone]{Corollary}
\newtheorem{kconv}[Emonotone]{Corollary}
\newtheorem{kslineq}[Emonotone]{Lemma}
\newtheorem{kssbound}[Emonotone]{Corollary}
\newtheorem{hbound}[Emonotone]{Corollary}
\newtheorem{L2bounds}[Emonotone]{Corollary}
\newtheorem{L2decay}[Emonotone]{Corollary}

As we remarked in the introduction, it is not clear that the quantity $E$, defined by \eqref{E:energy}, is necessarily monotone under the flow \eqref{r1}.  However, our first result of this section shows that in fact $E$ is monotone under the flow \eqref{r1} if it is initially sufficiently small.

\begin{Emonotone} \label{T:Emonotone}
There exists a constant $\varepsilon>0$, depending only on $\omega$ and $L_0$ and a constant $c>0$ depending only on $\omega$ such that, if 
\begin{equation} \label{E:smallE}
  E\left[ \gamma_0\right] < \varepsilon
\end{equation}
then, under the flow \eqref{r1} we have
$$\frac{d}{dt} E\left[ \gamma_t\right] \leq -\frac{c}{2} \int k_{s^4}^2 ds \mbox{.}$$
\end{Emonotone}

\noindent \textbf{Proof:} First note that using Lemma \ref{d12}, \eqref{r1} and \eqref{h2} we have 
\begin{multline} \label{E:dE}
  \frac{d}{dt} E = - \int \left( k_{s^4} + k^2 k_{ss} - \frac{1}{2} k\, k_s^2 \right) \left( G + h \right) ds \\= - \int G^2 ds - h\left( t\right) \int G\, ds
   = - \int G^2 ds + \frac{1}{2\pi \omega} \int G\, ds \int k\, G\, ds \mbox{.}
   \end{multline}
   
  We estimate using integration by parts and Lemma \ref{T:ksup}
\begin{multline*}
  \left| \int G\, ds \right| = \left| \int \left(  k_{s^4} + k^2 k_{ss} - \frac{1}{2} k\, k_s^2 \right) ds \right|\\
  \leq \frac{5}{2} \left\| k \right\|_\infty \int k_s^2 ds \leq 5 E\left( L_0^\frac{1}{2} E + \frac{2\pi\omega}{L_0} \right) =: F_1\left( E, \omega, L_0 \right) \mbox{.}
  \end{multline*}
  Therefore via the Peter-Paul inequality
  \begin{align*}
  \frac{1}{2 \pi \omega} \left| \int G\, ds \right| \left| \int k\, G\, ds \right| 
  &\leq \frac{F_1\left( E, \omega, L_0 \right)}{2 \pi \omega} \left| \int k\, G\, ds \right| \\
  &\leq \frac{1}{2} \int G^2 ds + \frac{F_1^2\left( E, \omega, L_0 \right)}{8 \pi^2 \omega^2} \int k^2 ds\\
  &\leq  \frac{1}{2} \int G^2 ds +  \frac{F_1^2\left( E, \omega, L_0 \right)}{8 \pi^2 \omega^2} \left\| k \right\|_\infty \cdot 2\pi \omega \\
  &=  \frac{1}{2} \int G^2 ds + \frac{F_1^2\left( E, \omega, L_0 \right)}{4 \pi \omega} \left( L_0^{\frac{1}{2}} E + \frac{2\pi \omega}{L_0} \right) 
  \end{align*}
  and from \eqref{E:dE} we obtain
  $$ \frac{d}{dt} E \leq -\frac{1}{2} \int G^2 ds + \frac{F_1^2\left( E, \omega, L_0 \right)}{4 \pi \omega} \left( L_0^{\frac{1}{2}} E + \frac{2\pi \omega}{L_0} \right) \mbox{.}$$
  Now
  $$\frac{F_1^2\left( E, \omega, L_0 \right)}{4 \pi \omega} \left( L_0^{\frac{1}{2}} E + \frac{2\pi \omega}{L_0} \right) \leq C_0\left( L_0, \omega \right) E^2\left( 1 + E^3 \right)$$
  so using also Corollary \ref{T:G2new} we have
  $$ \frac{d}{dt} E \leq -c \int k_{s^4}^2 ds + c_2\left( L_0, \omega \right) E^2\left( 1 + E^3 \right) \mbox{.}$$
  From Lemma \ref{T:PSW} we have
  $$E \leq \frac{1}{2} \left( \frac{L_0}{2\pi \omega}\right)^6 \int k_{s^4}^2 ds$$
  which implies
  $$\frac{d}{dt} E \leq -\left[ c - C_3\left( L_0, \omega \right) E\left( 1 + E^3\right) \right] \int k_{s^4}^2 ds \mbox{.}$$
  Thus for $E_0$ small enough, depending only on $L_0$ and $\omega$, we conclude the result.\hspace*{\fill}$\Box$\\

Using now Lemma \ref{T:PSW} similarly as in the previous proof, we obtain exponential decay of $E$ under the flow \eqref{r1}.

\begin{Edecay} \label{T:Edecay}
For $E\left[ \gamma_0\right]$ satisfying \eqref{E:smallE}, there exists a $\delta>0$, depending on $L_0$ and $\omega$ such that, under the flow \eqref{r1},
$$E\left[ \gamma_t\right] \leq E\left[ \gamma_0\right] e^{-\delta\, t} \mbox{.}$$
\end{Edecay}

We can also now establish an $L^2$ bound on the curvature of $\gamma_t$ under \eqref{r1} with small initial $E$.

\begin{kbound} \label{T:kbound}
For $E\left[ \gamma_0\right]$ satisfying \eqref{E:smallE}, there exists a $C_0>0$, depending only on $L_0$ and $\omega$ such that, under the flow \eqref{r1},
\begin{equation} \label{E:kbound}
  \left\| k \right\|_\infty \leq 2 L_0^{\frac{1}{2}} E\left[ \gamma_0\right] + \frac{2\pi \omega}{L_0} =: C_0\left( L_0, \omega, E\left[ \gamma_0\right] \right)\mbox{.}
\end{equation}
Since $L$ is constant, $\int k^2 ds$ is also uniformly bounded under \eqref{r1}.
\end{kbound}

\noindent \textbf{Proof:} In view of Corollary \ref{T:Edecay}, this follows immediately from Lemma \ref{T:ksup}.\hspace*{\fill}$\Box$\\

Moreover, using Lemma \ref{T:PSW}, since $L$ is fixed we immediately have from Corollary \ref{T:Edecay} pointwise exponential decay of the curvature to that of a circle of circumference length $L$.

\begin{kconv} \label{T:kconv}
For $E\left[ \gamma_0\right]$ satisfying \eqref{E:smallE}, there exist positive constants $C$ and $\delta$, depending on $L_0$ and $\omega$ such that, under the flow \eqref{r1},
$$\left\| k - \overline{k} \right\|_{\infty}^2 \leq C\, e^{-\delta\, t} \mbox{.}$$
\end{kconv}

Further estimates are required to establish smooth convergence of solutions to a unique limiting circle.  We next bound all derivatives of curvature, under \eqref{r1} with small `energy', beginning with the following evolution inequality.

\begin{kslineq} \label{T:kslineq}
For any $\eta>0$, we have under \eqref{r1},
$$\frac{d}{dt} \int k_{s^\ell}^2 ds \leq -\left( 2 - \eta \right) \int k_{s^{\ell+3}}^2 ds + c\left( \eta\right) \left( \int k^2 ds\right)^{2\ell + 7} - h\left( t \right) \int k\, k_{s^\ell}^2 ds \mbox{.}$$
\end{kslineq}

\noindent \textbf{Proof:} We estimate the $P$ terms in Lemma \ref{T:kslevln} using Lemma \ref{i13}.  For the first, by integrating by parts as necessary we can ensure no higher derivative than $k_{s^{\ell+2}}$ appears and thus 
$$2\int k_{s^\ell} P_3^{\ell+4}\left( k\right) \leq \frac{\eta}{2} \int k_{s^{\ell+3}}^2 ds + c\left( \eta\right) \left( \int k^2 ds\right)^{2\ell + 7} \mbox{.}$$
The second estimate follows similarly; resubstituting then gives the desired inequality.\hspace*{\fill}$\Box$\\

Since the highest derivative occurring in $h\left( t\right)$ is $k_{ss}$ (see \eqref{h2}), the crucial step is to bound $\int k_{ss}^2 ds$ under the flow.  This is what we establish next.

\begin{kssbound} \label{T:kssbound}
  For $E\left[ \gamma_0\right]$ sufficiently small, there exists a constant $C_2>0$ such that, under the flow \eqref{r1},
  $$\int k_{ss}^2 ds \leq C_2 \mbox{.}$$
  \end{kssbound}
  
  \noindent \textbf{Proof:} From Lemma \ref{T:kslineq} and \eqref{h2} we have
\begin{multline*}
  \frac{d}{dt} \int k_{ss}^2 ds \leq -\left( 2 - \eta \right) \int k_{s^{5}}^2 ds + c\left( \eta\right) \left( \int k^2 ds\right)^{11} + \frac{1}{2\pi \omega} \int k_{ss}^2 ds \int k\, k_{ss}^2 ds\\
   - \frac{7}{4\pi \omega} \int k^2 k_s^2ds \int k\, k_{ss}^2 ds \mbox{.}
   \end{multline*}
  Using now Lemma \ref{T:kbound} together with integration by parts and the H\"{o}lder inequality, we estimate
  $$ \frac{1}{2\pi \omega} \int k_{ss}^2 ds \int k\, k_{ss}^2 ds \leq \frac{C_0}{2\pi \omega} \left( \int k_{ss}^2 ds\right)^2 =\frac{C_0}{2\pi \omega} \left( -\int k_s k_{s^3} ds\right)^2 \leq \frac{C_0}{\pi \omega} E\left[ \gamma_t\right] \int k_{s^3}^2 ds$$
  and 
  $$- \frac{7}{4\pi \omega} \int k^2 k_s^2ds \int k\, k_{ss}^2 ds \leq \frac{7 C_0^3}{2\pi} E\left[ \gamma_t\right] \int k_{ss}^2 ds \mbox{.}$$
  Hence
  $$\frac{d}{dt} \int k_{ss}^2 ds \leq -\left( 2 - \eta \right) \int k_{s^{5}}^2 ds + c\left( \eta\right) \left( \int k^2 ds\right)^{11} + \frac{C_0}{\pi \omega} E\left[ \gamma_t\right] \int k_{s^3}^2 ds +\frac{7 C_0^3}{2\pi} E\left[ \gamma_t\right] \int k_{ss}^2 ds \mbox{.}$$
  Finally using Lemma \ref{T:PSW} and Proposition \ref{T:Emonotone} we obtain
  $$\frac{d}{dt} \int k_{ss}^2 ds \leq -\left( 2 - \eta - \frac{C_0 E\left[ \gamma_0\right] L_0^4}{16\pi^5} - \frac{7 C_0^3 E\left[ \gamma_0 \right] L_0^6}{128 \pi^7} \right) \int k_{s^5}^2 ds + c\left( \eta\right) \left( \int k^2 ds\right)^{11} \mbox{.}$$
  Assuming $E\left[ \gamma_0\right]$ is small enough, it follows by choosing suitably small $\eta$ that the leading coefficient is negative, and in view of Corollary \ref{T:kbound}, the result follows. \hspace*{\fill}$\Box$\\

\begin{hbound} \label{T:hbound}
For $E\left[ \gamma_0\right]$ sufficiently small as in Corollary \ref{T:kssbound}, there is a corresponding constant such that, under the flow \eqref{r1},
$$\left| h\left( t\right) \right| \leq C_h \mbox{.}$$
\end{hbound}

\noindent \textbf{Proof:} We estimate from \eqref{h2}
$$\left| h\left( t\right) \right| \leq \frac{1}{2\pi} \int k_{ss}^2 ds + \frac{7}{4\pi} \int k^2 k_s^2 ds$$
from which the result follows using Corollaries \ref{T:kssbound}, \ref{T:kbound} and Proposition \ref{T:Emonotone}.\hspace*{\fill}$\Box$\\

We may now obtain bounds and exponential decay in $L^2$ and $L^\infty$ for all curvature derivatives by standard arguments.

\begin{L2bounds} \label{T:L2bounds}
For $E\left[ \gamma_0\right]$ sufficiently small as in Corollary \ref{T:kssbound}, for each $\ell \in \mathbb{N} \cup \left\{ 0 \right\}$ there is a corresponding positive constant $C_\ell$ such that, under the flow \eqref{r1},
$$\int k_{s^\ell}^2 ds \leq C_\ell$$
as long as the solution exists.
\end{L2bounds}

\noindent \textbf{Proof:} We have already proven the cases $\ell = 0, 1, 2$.  For general $\ell$, using Corollary \ref{T:hbound} in Lemma \ref{T:kslineq} we have under \eqref{r1} for any $\eta>0$, 
$$\frac{d}{dt} \int k_{s^\ell}^2 ds \leq -\left( 2 - \eta \right) \int k_{s^{\ell+3}}^2 ds + c\left( \eta\right) \left( \int k^2 ds\right)^{2\ell + 7} +C_h  \int P_{3}^{2\ell}\left( k\right) ds \mbox{.}$$
For each $\ell$, with small enough $\eta$, the result now follows using Lemmas \ref{i13} and \ref{T:PSW}.\hspace*{\fill}$\Box$

\begin{L2decay}
For $E\left[ \gamma_0\right]$ sufficiently small as in Corollary \ref{T:kssbound}, for all $\ell \in \mathbb{N}$ there are corresponding positive constants $\tilde C_\ell, \delta$ such that, under the flow \eqref{r1},
$$\int k_{s^\ell}^2 ds \leq \tilde C_\ell e^{-\delta \, t}$$
as long as the solution exists.
\end{L2decay}

\noindent \textbf{Proof:} For $\ell=1$ this is Corollary \ref{T:Edecay}.  For $\ell>1$ the result follows inductively integration by parts, the H\"{o}lder inequality, Corollary \ref{T:Edecay}, Proposition \ref{T:kbound} and Corollary \ref{T:L2bounds}.\hspace*{\fill}$\Box$\\
  
  Arguments as in \cite[Section 7]{AMWW} now complete the proof of Theorem
\ref{T:main}.  In particular, exponential decay of curvature derivatives in
$L^\infty$ facilitates boundedness of $\gamma$ and all its derivatives, giving
boundedness in $C^\infty$ of $\gamma\left( \cdot, t_j \right)$ along any
sequence $t_j\rightarrow \infty$.  Exponential decay of $E\left[ \gamma\left(
\cdot, t\right)\right]$ implies convergence of $\gamma\left( \cdot, t_{j_k}
\right)$ to a circle of radius $\frac{L_0}{2\pi\omega}$, in the $C^\infty$-topology,
for a subsequence $\left\{ t_{k_j} \right\}$.  Full convergence of the flow to
the $\omega$-circle, modulo reparametrisation, follows via the result \cite[Appendix
A]{AMWW}.\hspace*{\fill}$\Box$

\end{document}